\documentclass[12pt,draftcls,onecolumn]{IEEEtran}
\newtheorem{theorem}{Theorem}

\newtheorem{lemma}[theorem]{Lemma}

\begin{document}
%
\title{On the NP-Hardness of Checking
Matrix Polytope Stability and Continuous-Time Switching Stability}
%
%
%


\author{Leonid Gurvits and Alex Olshevsky
\thanks{L. Gurvits is with Los Alamos National Laboratory, Los Alamos, NM, 87545.
        {\tt\small gurvits@lanl.gov}}%
\thanks{A. Olshevsky is with the Laboratory for Information and Decision Systems, Massachusetts Institute of Technology, Cambridge,
MA, 02139. {\tt\small alex\_o@mit.edu}}
\thanks{A preliminary version of this paper has appeared in the \em{Proceedings of the European Control Conference, Kos, Greece, 2007.}}%
\thanks{Manuscript received ?????; revised ?????. }}

%
%

\markboth{IEEE Transactions on Automatic Control,~Vol.~????, No.~????, ?????}%
{Shell \MakeLowercase{\textit{et al.}}: IEEE Transactions on
Automatic Control,~Vol.~????, No.~????, ?????}
%



\maketitle

\begin{abstract}
Motivated by questions in robust control and switched linear
dynamical systems, we consider the problem checking whether all
convex combinations of $k$ matrices in $R^{n \times n}$ are stable.
In particular, we are interested whether there exist algorithms
which can solve this problem in time polynomial in $n$ and $k$. We
show that if $k= \lceil n^d \rceil$ for any fixed real $d>0$, then
the problem is NP-hard, meaning that no polynomial-time algorithm in
$n$ exists provided that $P \neq NP$, a widely believed conjecture
in computer science. On the other hand, when $k$ is a constant
independent of $n$, then it is known that the problem may be solved
in polynomial time in $n$.  Using these results and the method of
measurable switching rules, we prove our main statement: verifying
the absolute asymptotic stability of a continuous-time switched
linear system with more than $n^d$ matrices $A_i \in R^{n \times n}$
satisfying $0 \succeq A_i + A_i^{T}$ is NP-hard.
\end{abstract}


%

\section{Introduction}
Let $\{A_1, \ldots, A_k\}$ be a finite set of real $n \times n$
matrices. We define the corresponding matrix polytope ${\cal A}$ as
those matrices which can be written as $\sum_{i=1}^k \alpha_i A_i$
for some nonnegative real numbers $\alpha_1,\ldots, \alpha_k$ adding
up to $1$. We are concerned with the following decision problem,
which we will call the $(k,n)$-POLYTOPE-STABILITY problem: given $k$
rational $n \times n$ matrices $\{A_1,...,A_k\}$ to decide whether
every matrix in $\cal{A}$ is Hurwitz stable, i.e. has eigenvalues
with negative real parts.

The $(k,n)$-POLYTOPE-STABILITY problem, and some of its natural
generalizations, have been considered before in the control theory
literature (see \cite{B85}, \cite{FB88}, \cite{BFS88}, \cite{S90},
\cite{W91}, \cite{CL93}, \cite{B95}, \cite{M99}, \cite{HAPL04},
\cite{B04}, \cite{TB04}, \cite{CGTV05}, \cite{C05},\cite{OP05},
\cite{S05}). Notable results include the solution of the $k=2$ case
in \cite{B85},\cite{FB88}, \cite{S90}, also known as the stability
testing of affine representations; a Lyapunov-search algorithm in
\cite{CL93}; and a recent approach based on LMI relaxations
\cite{HAPL04}, \cite{CGTV05}, \cite{TB04}.

Our first result is that when there are $\lceil n^d \rceil$ different $n \times n$
rational matrices $A_i$, where $d$ is any positive real number (in
the above notation, the $(\lceil n^d \rceil, n)$-POLYTOPE-STABILITY
problem) the problem of deciding whether there exists an unstable
matrix in $\cal{A}$ is NP-hard.

On the other hand,  in many circumstances where the
$(k,n)$-POLYTOPE-STABILITY problem appears, $k$ is constant that
does not vary with $n$. Consider, for example, testing the stability
of an
               $n \times n$ interval matrix where all but $r$ of the entries are
               known precisely. It is easy to see that this is a special case of
               the polytope stability problem with $k=2^r$. Moreover, if $r$ is
               fixed, but $n$ is allowed to vary, we end up with a polytope
               stability problem with fixed $k$.

We note that the $(k,n)$-POLYTOPE-STABILITY problem can be solved in
$n^{O(k)}$ elementary operations by reducing the problem to deciding
whether a certain multivariate polynomial does not have real roots,
which can in turn be solved using quantifier elimination algorithms
- see \cite{GS07} for a detailed writeup (the reduction to a root
localization problem for a multivariate polynomial was first done in
\cite{HAPL04}). If the number $k$ is fixed, this gives a
deterministic polynomial-time algorithm.



The second and main subject of our paper is the $(k,n)$ - Continuous Time
Absolute Switching Stability Problem ($(k,n)-CTASS$ problem): given
$k$ rational matrices $A_1,\ldots,A_k$ in $R^{n \times n}$ to decide
whether there exists a norm $||\cdot||$ in $R^n$ and $a
> 0$ such that the induced operator norms satisfy the inequalities:
\begin{equation} \label{lyapun}
||exp(A_i t)|| \leq e^{-at} :  a > 0, ~t \geq 0 , ~1 \leq i \leq k.
\end{equation}



We consider matrices that satisfy the further restriction of being
non-strict contractions in the two-norm: \[ ||exp(A_i t)||_{2} \leq
1 , t \geq 0, 1 \leq i \leq k, \] which is equivalent to requiring
\[ 0 \succeq A_i + A_i^{*}. \]

We show that, even in this case, checking the condition of Eq.
(\ref{lyapun}) with more than $n^d$ matrices $A_i \in R^{n \times n}$, where $d$ is
any positive real number,
is also NP-hard. As far as we know, this is the first hardness
result for continuous-time switched linear systems.\\
{We note that it is well known that the stability of the polytope
${\cal A}$ is necessary, but not in general sufficient, for the
absolute switching stability condition of Eq. (\ref{lyapun}) to
hold. Luckily, these two conditions are equivalent for the
``gadgets" used in our proof of NP-HARDNESS of polytope stability.}

We stress that our result says nothing about the solvability of
concrete, finite-size instances for the problem - for instance, this
paper has no new implications
  for testing the stability of all convex combinations of $3$ matrices in $R^{8 \times 8}$. Rather our result
implies that if $P \neq NP$ then any algorithm for solving this
problem with $n^d$ matrices $n \times n$ matrices, for any $d > 0$,
will have a worst-case operations count which grows faster than any
polynomial in $n$.

This provides an explanation why many approaches to this problem
fail. Despite an extensive literature devoted to this problem cited
above, no polynomial time algorithms are known. See, in particular,
\cite{BFS88}, for proofs that various intuitive approaches fail. Our
result implies that in fact any polynomial time algorithm for this
problem would immediately disprove the $P \neq NP$ conjecture.

\section{NP-Hardness of Stability Testing of Matrix Polytopes \label{hardnesssection}}

In this section, we consider the computational complexity of
deciding whether every matrix in the set $\cal{A}$ (defined by Eq.
(\ref{polydef})) is stable. We will show that this problem is
NP-hard through a reduction from the maximum clique problem, which
is known to be NP-complete \cite{K72}. The details of this reduction
are described below.

\begin{theorem} $(n,2n+2)$-POLYTOPE-STABILITY is NP-hard. \label{hardnesstheorem} \end{theorem}

Notice that the well-known interval stability problem\footnote{The
interval stability problem is to determine, given numbers
$\{\underline{a}_{ij},\bar{a}_{ij}\}_{i,j=1 \ldots m}$, whether
every matrix $A$ satisfying $A_{ij} \in
[\underline{a}_{ij},\bar{a}_{ij}]$
 is stable.} corresponds to
$(k,m)$-POLYTOPE-STABILITY with $k$ exponential in $m$. In other
words NP-hardness of interval stability, shown in \cite{N93} and
\cite{PR93}, does not imply NP-hardness of
$(n,2n+2)$-POLYTOPE-STABILITY.

We first give a series of definition and lemmas before proving
Theorem \ref{hardnesstheorem}. Given $k$ rational matrices
$A_1,\ldots, A_k$ in $R^{n \times n}$, we will refer to the problem
of deciding whether there exists a singular matrix in the associated
polytope $\cal{A}$ as the $(k,n)$-POLYTOPE-NONSINGULARITY problem.

\begin{lemma} There is a polynomial-time reduction from the
$(k,n)$-POLYTOPE-NONSINGULARITY problem to the
$(k,2n)$-POLYTOPE-STABILITY problem \label{stable-singular}
\end{lemma}

\begin{proof} Given a square matrix $A \in \mathbf{R}^{n \times n}$ define
$B$ as
\[ B = \left(
                     \begin{array}{cc}
                       0_{n \times n} & A^T \\
                       -A & -I_{n \times n} \\
                     \end{array}
                   \right). \] We claim $B$ is Hurwitz if and only if $A$ is nonsingular.
This follows because spectrum of $B$ can be written as $\sigma(B) =
\cup_{1 \leq i \leq n} \{ \frac{1}{2} (-1 \pm \sqrt{1 -
(s_{i}(A))^{2}}) \}$, where $s_{1}(A) \leq ... \leq s_{n}(A)$ are
singular values of the matrix $A$.

 Suppose we are given $k$ $n \times n$ matrices $A_i$, and we want to decide
 whether the set $\cal{A}$ defined by Eq. (\ref{polydef}) contains a nonsingular matrix.  Define \[ B_i = \left(
                     \begin{array}{cc}
                       0_{n \times n} & A_i^T \\
                       -A_i & -I_{n \times n} \\
                     \end{array}
                   \right). \] Since
\[ \sum_i \alpha_i B_i = \left( \begin{array}{cc}
                       0_{n \times n} & (\sum_i \alpha_i A_i)^T \\
                       -\sum_i \alpha_i A_i & -I_{n \times n} \\
                     \end{array}
                   \right). \] when $\sum_i \alpha_i = 1$, it follows
                   by the previous item that testing
                   POLYTOPE-NONSINGULARITY with the set $\cal{A}$
is the same as testing POLYTOPE-STABILITY on the set
\[ {\cal B} = \{ B ~|~ B = \sum_i \alpha_i B_i, \sum_i \alpha_i = 1,
\alpha_i \geq 0 \mbox{ for all } i \}. \] However, note that the
construction has doubled the dimension, since the matrices $B_i$
belong to $R^{2n \times 2n}$. This concludes the proof that
$(k,n)$-POLYTOPE-NONSINGULARITY can be reduced to
$(k,2n)$-POLYTOPE-STABILITY. \end{proof}

Consider the problem of deciding whether there
                   exists a nonnegative vector $p$ in $R^n$ whose components sum to $1$
                   such that $p^T M p = 1$ for an arbitrary invertible
                   matrix $M$. We will consider $M^{-1}$ to be the input to this problem. We will refer to the problem as the
                   $n$-QUADRATIC-THRESHOLD problem.

\begin{lemma} There is a polynomial-time reduction from the $n$-QUADRATIC-THRESHOLD problem
to the $(n,n+1)$-POLYTOPE-NONSINGULARITY problem.
\label{threshold-singular}\end{lemma}

\begin{proof}

Define \[ X_i^{(n)} = \left(
                  \begin{array}{cc}
                    M^{-1} & e_i \\
                    e_i^T &  1\\
                  \end{array}
                \right),
 \] where $e_i$ is the column vector with $1$ in the
i'th entry and zeros elsewhere. Define
  $ {\cal S}_n  = \{p \in R^n ~|~ \sum_i p_i = 1, ~p_i \geq 0 ~\forall i
  \}$ and let

\[ {\cal X} = \{ X ~|~ X = \sum_i p_i X_i, ~p \in {\cal S}_n \}. \] In other words, $\cal{X}$ is
the set of matrices of the form
\[ \left(
     \begin{array}{cc}
       M^{-1} & p^T \\
       p &  1\\
     \end{array}
   \right), \] with $p \in {\cal S}_n$. By the Schur complement formula such a matrix is singular if and only if $p^T
   M p = 1$. Thus given an invertible matrix $M$, we can solve the $n$-QUADRATIC-THRESHOLD problem by solving an instance of the $n+1$-POLYTOPE-NONSINGULARITY
   problem with the polytope $\cal{X}$. \end{proof}

   Let $G=(V,E)$ be an undirected graph. A subset of the vertices $C$ is called a clique
   if $v_1,v_2 \in C$ implies $(v_1,v_2) \in E$. The $n$-MAX-CLIQUE problem is the problem of determining the
   size of the largest clique, denoted by $\omega(G)$, in an undirected graph $G$ on $n$ vertices.

\begin{lemma} There is a polynomial-time reduction from the $n$-MAX-CLIQUE  problem to the
$n$-QUADRATIC-THRESHOLD problem \label{clique-threshold}\end{lemma}

\begin{proof}

\noindent {\bf 1.} It is known that \cite{MS65}:

   \begin{equation} \label{cliqueeq}  1-\frac{1}{\omega(G)} = \max_{p \in {\cal S}_n} p^T M p.
   \end{equation}   where $M$ is the adjacency matrix of the graph $G$.

\noindent {\bf 2.} Because the QUADRATIC-THRESHOLD problem is
defined only for nonsingular matrices, for our reduction to work we
will need to modify $M$ to insure its nonsingularity. To this end,
we consider the matrices $M_i = M + \frac{1}{n^2+i} I$ for
   $i=1,\ldots,n+1$. At least one $M_i$ must be nonsingular, because
   $M$ cannot have $n+1$ eigenvalues. We find a nonsingular $M_i$ (this can be done
   in polynomial time with Gaussian elimination for each $i=1,\ldots,n+1$). Let us denote this nonsingular $M_i$ by $M_{i^*}$.\\

   In the proof below, we will threshold the form $p^T M_{i^*} p$;
   recall from Lemma \ref{threshold-singular} that this requires
   the computation of $M_{i^*}^{-1}$. This involves a polynomial number of computations in $n$, and
   the bit-sizes remain polynomial as well; for a proof see Corollary 3.2a of \cite{S86}.

\noindent {\bf 3.}   We have that for $p \in {\cal S}_n$,
\[ p^T M_{i^*} p  =  p^T M p + \frac{1}{n^2+i^*} \sum_i
   p_i^2 \]Because $p_i \in [0,1]$ and $\sum_i p_i = 1$ imply that $\sum_i p_i^2 \leq 1$, we have:

   \begin{equation} p^T M p  \leq  p^T M_{i^*} p \leq p^T M p +
   \frac{1}{n^2+i^*} \label{bounds}
   \end{equation}

   It follows that \begin{equation} \label{modthresh} 1-\frac{1}{\omega(G)} \leq \max_{p \in {\cal S}_n} p^T M_{i^*} p
   \leq 1 - \frac{1}{\omega(G)} +\frac{1}{n^2+i^*} \end{equation}

\noindent {\bf 4.}     Because the optimal solution of Eq.
(\ref{cliqueeq}) is $1-1/\omega(G)$, and $\omega(G)$ is an integer
between $1$ and $n$,
   this optimal solution must be in the
   set $S = \{0,1/2,2/3,3/4,\ldots,1-1/n\}$. Because the gap between the elements of $S$ is less than $1/n^2$, and consequently less than $1/(n^2+i^*)$, it follows from Eq. (\ref{modthresh})
   that the largest element of $S$ smaller than $\max_{p \in {\cal S}_n} p^T M_{i^*} p$ must be $1-\frac{1}{\omega(G)}$, the solution of the
   MAX-CLIQUE problem. Let this element be called $k^{*}$; then, it
   follows that \begin{equation} \label{qt}
  \frac{1}{k^*} p^T M_{i^*}
   p \geq 1. \end{equation} for some $p \in {\cal S}_n$, and $k^{*}$ is the largest element of $S$ with this property.

   For each $k \in S$, the existence of a $p \in {\cal S}_n$ satisfying Eq. (\ref{qt}) can, due to the invertibility of
   $M_i^{*}$, be decided by evaluating $p^T M_{i^*} p$ at an arbitrary $p \in {\cal S}_n$, followed up with a call to the QUADRATIC-THRESHOLD
   problem.  This is the reduction from MAX-CLIQUE to QUADRATIC-THRESHOLD.  \end{proof} \vspace{0.1cm}

\begin{proof}[Proof of Theorem \ref{hardnesstheorem}]: Lemmas
\ref{stable-singular}, \ref{threshold-singular},
\ref{clique-threshold} provide a reduction from the MAX-CLIQUE
problem to the POLYTOPE-STABILITY problem. The size of the problem
goes from $n$ in the QUADRATIC-THRESHOLD problem to $(n,n+1)$ after
the reduction to POLYTOPE-NONSINGULARITY; and from $(n,n+1)$ to
$(n,2(n+1))$ in the reduction from POLYTOPE-NONSINGULARITY to
POLYTOPE-STABILITY. Since MAX-CLIQUE is known to be NP-complete
\cite{K72}, it follows that
   $(n,2n+2)$-POLYTOPE-STABILITY is NP-hard. \end{proof}

\noindent {\bf Remark:} Note that the matrices $A_i$ created in our
reduction have entries whose bit-size are polynomial in $n$.

\noindent {\bf Remark:}
\begin{enumerate}
\item
We note that our results easily imply the
NP-HARDNESS of POLYTOPE-STABILITY with $\lceil n^d \rceil $ extreme points, for any
real $d>0$ (here $\lceil x \rceil$ refers to the smallest integer which is at least
$x$).  Indeed, the $\lceil n^d \rceil$-MAX-CLIQUE
problem remains NP-COMPLETE,
and therefore the $(\lceil n^d \rceil,2 \lceil n^d \rceil+2)$-POLYTOPE-STABILITY remains
NP-Hard. However, we clearly do not make the problem any easier by
increasing the dimension, so that the $(\lceil n^d \rceil,2n+2)$ problem is
NP-Hard.
\item We have remarked that the $(k,n)$-POLYTOPE STABILITY problem may be solved
in polynomial time if $k$ is upper bounded by a constant. We may ask
about the reverse question: what happens when $n$ is upper bounded
by a constant? Is the problem solvable in polynomial time as a
function of $k$?

The answer is yes. Caratheodory's theorem implies that any matrix in ${\cal A}$ may be
expressed as a convex combination of $n^2+1$ matrices. Thus we reduce the problem of
checking $(k,n)$ polytope stability to the problem of checking ${k \choose n^2+1}$
different $(n^2+1,n)$ polytope stability problems. When $n$ is upper bounded,
checking $(n^2+1,n)$-POLYTOPE-
STABILITY (say by computing the determinant explicitly and using quantifier elimination)
takes a constant number of operations, so the number
of operations grows as ${k \choose n^2+1}$,
which, when $n$ is upper bounded, is polynomial in $k$.

\item
Essentialy, our construction boils down to the next determinantal representation:
\begin{equation} \label{val}
Q(p_1,...,p_n) = \det(A_0 + \sum_{1 \leq i \leq n} p_i A_i)
\end{equation}
Such representations exists for any polynomial $Q(p_1,...,p_n)$ \cite{val}.
In our case\\
$Q(p_1,...,p_n) = <Mp,p> -1$.
\end{enumerate}
 \section{NP-HARDNESS of Checking Continuous-Time Absolute Switching Stability}
In this section, we will show that checking the absolute switching
stability of a class of continuous-time linear switched systems is
NP-hard. Recall that the $(k,n)$ continuous time switching stability
problem is: given $k$ rational matrices $A_1,\ldots,A_k$ in $R^{n
\times n}$, to decide whether there exists a norm $||\cdot||$ in
$R^n$ and $a
> 0$ such that the induced operator norms satisfy the inequalities:
\begin{equation} \label{normeq}
||exp(A_i t)|| \leq e^{-at} : t \geq 0 , 1 \leq i \leq k.
\end{equation}

We consider a subcase of the problem where the matrices $A_i$
satisfy the (nonstrict) Lyapunov inequalities $0 \succeq A_i
+A_i^{T}, 1 \leq i \leq k$. In assuming this, we only make the
problem easier, as we assume that the matrices $A_i$ have a nice
geometrical structure; indeed, the previous requirement corresponds
to requiring that solutions of the equation $\dot{x}(t)=A_i x(t), i
= 1,\ldots,k$ with measurable switching rules, are nonincreasing in
the $2$-norm.

Nevertheless, we will show that testing continuous time stability is
NP-hard already in this case.  The following lemma -   together with
Theorem \ref{hardnesstheorem} - proves this. As far as we know it is
the first hardness result in the area of continuous time absolute
switching stability .

\begin{lemma} Consider the following $2n \times 2n$ matrices $B_i , 1 \leq i \leq k < \infty$:

\[ B_i = \left(
                     \begin{array}{cc}
                       0_{n \times n} & A_i^T \\
                       -A_i & -I_{n \times n} \\
                     \end{array}
                   \right)\]
Then there exists a norm $||\cdot||$ in $R^{2n}$ and $a > 0$ such
that the induced operator norms satisfy the following inequalities:
\begin{equation}
||exp(B_i t)|| \leq e^{-at} : t \geq 0 , 1 \leq i \leq k ,
\label{decayeq} \end{equation} if and only if all matrices in the
convex hull $\cal{A}$ are nonsingular. \label{switch}\end{lemma}
\bigskip

Before proving the lemma, we need the following auxiliary claim.
Consider the following family of "differential"
equations\footnote{Strictly speaking, they ought to be viewed as
integral equations \[ x(t) = x(0) + \int_{0}^{t}([\sum_{i=1}^k
p_i(\tau) B_i] x(\tau))d \tau.\]}:
$$
\dot{x}(t) = [\sum_{i=1}^k p_i(t) B_i] x(t),
$$ with initial condition satisfying $||x(0)||_2=1$. Here $\left(
                                                              \begin{array}{ccc}
                                                                p_1(t) & \ldots & p_k(t) \\
                                                              \end{array}
                                                            \right)$
is a Lebesgue-measurable vector function whose range is a subset of
${\cal S}_k$. Since the $\sum_i p_i B_i$ lies in a bounded set of
matrices, the above equation has a unique Lipschitz solution. Since
$0 \succeq B_i + B_i^{T}$, it follows that $||x(t)||_{2} \leq 1$ for
$t \geq 0$.

\noindent {\bf Claim:} Assume that there is no induced norm
satisfying Eq. (\ref{normeq}). Then, there exists a measurable
vector function $\left(
                                                              \begin{array}{ccc}
                                                                p_1(t) & \ldots & p_k(t) \\
                                                              \end{array}
                                                            \right)$ whose range is a subset of ${\cal
                                                            S}_k$,
and a vector $x(0)$ with $||x(0)||_{2} = 1$ such that the solution
of
\begin{equation} \label{curve}
\dot{x}(t) = [\sum_{i=1}^k p_i(t) B_i] x(t),
\end{equation}
with initial condition $x(0)$ satisfies $||x(1)||_{2} = 1$.

\noindent {\bf Proof of claim:} We will prove the contrapositive of
the claim. It can be seen that the set of all possible values of
$x(t)$ at time $t=1$ produced by choices of $p(t),x(0)$ which
satisfy our assumptions is compact; see Theorem 4.7 in \cite{G02}
for details of the proof. So, suppose the conclusion is not true, by
compactness this means there exists $\delta
> 0$ such that for every $p,x(0)$ satisfying our assumptions, we
have that $||x(1)||_2 < 1-\delta$.  Thus there exists an
$\epsilon>0$ such that the system
\begin{equation}
\dot{x}(t) = [\sum_{i=1}^k p_i(t) (B_i+\epsilon I)] x(t),
\label{perturbedsystem} \end{equation}  has the same property( i.e.
there exists $\hat{\delta}$ with $||x(1)||_2 < 1-\hat{\delta}$ for
all suitable choices of $x(0),p(t)$).

We will define a norm such that $||e^{B_i t}|| \leq e^{-at}$ for
some $a>0$ and all $t \geq 0$, thus violating Eq. (\ref{decayeq}).
Define a norm on $R^n$ as follows. For $q \in R^n$,

\[ ||q||_{n} = \sup_{u \geq 1} \sup_{\mbox{ measurable } p(z), z \in [0,u] \mbox{ with range in } S_k}
||x(u)||_2 \] where $x(\cdot)$ is a solution to Eq.
(\ref{perturbedsystem}) with $x(0)=q$.

This norm induces a norm on $R^{n \times n}$:

\[ ||Q||_{n \times n} = \sup_{x \in R^n \\ ||x||_n = 1} ||Qx||_n
\]

However, \[ ||e^{(B_i + \epsilon I)t} x ||_{n} \leq ||x||_{n}
\] since premultiplication by $e^{B_i + \epsilon I}$ corresponds to
simply requiring that $p(z)=e_i$ for the first $t$ time units.
Therefore,

\[ ||e^{B_i t}||_{n \times n} \leq e^{- \epsilon t} \] which proves
the claim.

\begin{proof}[Proof of Lemma \ref{switch}:] We remark that the
argument is very similar to Theorem 4.7 in \cite{G02} and Corollary
2.8 in \cite{G03}.

First, we show the ``only if'' part, that is, assuming Eq.
(\ref{decayeq}), all the matrices in the convex hull $\cal{A}$ must
be nonsingular. Indeed, suppose not; suppose there exists a vector
$p$ such that $\sum_{i=1}^k p_i B_i$ is not stable, with
$\sum_{i=1}^k p_i = 1$ and all $p_i \geq 0$. Let $\lambda$ be an
eigenvalue of $\sum_i p_i B_i$ with nonnegative real part; then,
$e^{\lambda}$ will be an eigenvalue of $e^{\sum_i p_i B_i}$. Thus,
$e^{\sum_i p_i B_i}$ has an eigenvalue of magnitude at least $1$, so
that
\begin{equation} \label{eqbigger} ||e^{\sum_i p_i B_i}|| \geq 1,
\end{equation} where $||\cdot||$ is the same norm as in Eq. (\ref{decayeq}).
 However, by the Baker-Campbell-Hausdorff formula, for
all $\epsilon
> 0$, there exists $m$ large enough so that

\[ ||(e^{\sum_{i=1}^k p_i B_i})|| \leq (1+\epsilon) ||(\prod_{i=1}^k e^{(1/m) p_i
B_i})^m|| \]

Since $||e^{(1/m) p_i B_i}||^{m} \leq e^{-a p_i}$ by  Eq.
(\ref{decayeq}), we have that large enough $m$
\begin{equation} \label{eqsmaller}||e^{\sum_{i=1}^k p_i B_i}||
\leq (1+\epsilon) e^{-a \sum_i p_i} = (1+\epsilon) e^{-a} < 1,
\end{equation} where we pick $\epsilon$ small enough for the last inequality to hold. Equations
(\ref{eqbigger}) and (\ref{eqsmaller}) are in contradiction. We
conclude that all the matrices in $\cal{B}$ are indeed Hurwitz. By
Lemma \ref{stable-singular}, this implies all matrices in $\cal{A}$
are nonsingular. This proves the ``only if'' part.\\

Next, we show the ``if'' part. Let $x(0),p(t)$ be such that the
conclusion of the above claim
 is satisfied. The corresponding curve $x(t) : 0 \leq t \leq 1$, satisfying (\ref{curve}) is Lipschitz;
 $||x(t)||_{2} = 1$ for all $0 \leq t \leq 1$. It follows from the structure of matrices $B_i$ that
 if $C \in \cal{B}$ then for all non-zero vectors $x$ the inner product $\langle Cx,x \rangle \leq 0$ and $ \langle Cx,x \rangle = 0$ if and only if
 $x \in R^n \oplus 0_n$. Here $R^n \oplus 0_n$ is the subspace spanned by $\{e_1,\ldots,e_n\}$.
 If for some $\tau \in [0,1)$ the vector $x(\tau)$ does not belong to $R^n \oplus 0_n$ then it also holds
 by continuity of the curve $x(t)$ in a sufficiently small neighbourhood $[\tau, \tau + \epsilon]$.
 This implies that  $\langle x(\tau), C x(\tau + \delta) \rangle < -b < 0$ for some $b>0$ and all $0 \leq \delta \leq \epsilon$ and $C \in \cal{B}$.
 Define the vector $v(\epsilon) = \int_{\tau}^{\tau +\epsilon }([\sum_{i=1}^k p_i(s) B_i] x(s))d s$.
 It is clear that $||v(\epsilon)||_{2} \leq K \epsilon$ for some constant $K$ and that $\langle x(\tau), v(\epsilon) \rangle \leq -b \epsilon$.
 Therefore
 \begin{eqnarray*}
 ||x(\tau + \epsilon)||_2^{2} & = &  ||x(\tau)||_{2}^{2} + 2 <x(\tau), v(\epsilon)> + \\ &&  ||v(\epsilon)||_{2} \\ & \leq & 1 - 2b \epsilon + K^2 \epsilon^2.
 \end{eqnarray*}
 We get for small enough $\epsilon$ the inequality $||x(\tau + \epsilon)||_2 < 1$. We conclude
 that $x(t) \in R^n \oplus 0_n : 0 \leq t \leq 1$.

 This
 gives that $\int_{0}^{t} (\sum_{1 \leq i \leq k} p_i(\tau) A_i) x_1(\tau) d \tau = 0$ for all $0 \leq t \leq
 1$,
 where $x_1(\tau)$ is the vector formed by the first $n$ components of $x(\tau)$.
 As the Lebesgue measurable vector function on $[0,1]$ $\left(
                                                              \begin{array}{ccc}
                                                                p_1(t) & \ldots & p_k(t) \\
                                                              \end{array}
                                                            \right) \in {\cal S}_k $ is bounded
 thus  it follows that $(\sum_{1 \leq i \leq k} p_i(\tau) A_i) x_1(\tau) = 0$ up to measure zero.
 Since the last $n$ components of $x(\tau)$ are zero, and $||x(t)||_2=1$, we have  $||x_1(t)||_2=1 , 0 \leq t \leq 1$. Therefore there must exist a singular matrix in
 $\cal{A}$. \end{proof}

 \bigskip

\noindent {\bf Remark:}  In the discrete time case, it is known that
given two $n \times n$ rational matrices $A, B$ it is NP-hard to
check if there exists a norm $||.||$ in $R^n$ such that the induced
norms $||A||, ||B|| < 1$ \cite{TB97}. On the other hand, it had been
observed in \cite{G96}, that a slight modification of a construction
in \cite{TB97} gives a direct proof of the following statement:
given two $n \times n$ rational matrices $A, B$ with nonnegative
entries and $||A||_{l_1},||B||_{l_1} \leq 1$ it is NP-HARD to check
if there exists a norm $||.||$ in $R^n$ such that the induced norms
$||A||, ||B|| < 1$. However, the continuous-time counterpart of this
last problem is ``easy''(see Theorem 2.1 in \cite{G03}).Based on
this, it is unclear whether it is possible to modify the
constructions from \cite{TB97},\cite{G96} to handle the continious
time case.
\begin{small}

\end{small}

\begin{thebibliography}{99}



\bibitem{B85} S. Bialas, ``A Necessary and Sufficient Condition for
the Stability of Convex Combinations of Polynomials and Matrices,''
{\it Bulletin of the Polish Academy of Sciences. Technical
Sciences.} vol 33, pp. 473-480, 1985.

\bibitem{B95} B. Barmish, {\it New tools for robustness of linear systems,} New York:
MacMillan, 1995.

\bibitem{BPR96} S. Basu, R. Pollack, M.-F. Roy, ``On the
combinatorial and algebraic complexity of quantifier elimination,''
{\it Journal of the ACM,} vol. 43, p. 1002,  1996.

\bibitem{B04} P.-A. Bliman, ``A convex approach to robust stability
for linear systems with uncertain scalar parameters,'' {\it SIAM
Journal on Control and Optimization,} 42(6):2016–2042, 2004.

\bibitem{TB97} V. D. Blondel and J. N. Tsitsiklis, ``The Lyapunov exponent and joint spectral
radius of pairs of matrices are hard - when not impossible - to
compute and to approximate,'' {\it Mathematics of Control, Signals
and Systems,} Vol. 10, No. 1, 1997, pp. 31-40.


\bibitem{TB99} V. D. Blondel and J. N. Tsitsiklis, "A Survey of Computational Complexity Results in Systems and Control", {\it Automatica,} Vol. 36, No. 9, pp. 1249-1274, September 2000.

\bibitem{BFS88} R.B. Barmish, M. Fu, S. Saleh, ``Stability of a
polytope of matrices: counterexamples,'' {\it IEEE Transactions on
Automatic Control,} vol. 33, No. 6, pp. 569-572, 1988.

\bibitem{C05} G. Chesi, ``Establishing stability and instability of matrix
hypercubes,'' {\it Systems Control Letters,} 54(4):381–388, April
2005.

\bibitem{CGTV05} G. Chesi, A. Garulli, A. Tesi, and A. Vicino, ``Polynomially parameter-dependent
Lyapunov functions for robust stability of polytopic systems: An LMI
approach,'' {\it IEEE Transactions on Automatic Control,}
50(3):365–370, March 2005.

\bibitem{CL93} N. Cohen, I. Lewkowicz, ``A Necessary and Sufficient Criterion for Stability of Convex Sets of Matrices,''
{\it IEEE Transactions on Automatic Control}, vol. 38, No. 4, pp. 611-615, April 1993.

\bibitem{FB88} M. Fu, R. Barmish, ``Maximal unidirectional
perturbation bounds for the stability of polynomials and matrices,''
{\it Systems and Control Letters,}  vol. 11, pp. 173-179, 1988.

\bibitem{G96} L. Gurvits, "Stability of linear inclusions: Part 2", NEC Research report TR96-173 (1996).

\bibitem{G02} L. Gurvits, ``Stabilities and Controllabilities of Switched Systems (with Applications to the Quantum
Systems),'' {\it Proceeding of the 15th International Symposium on
the Mathematical Theory of Networks and Systems (MTNS '02), 2002}.

\bibitem{G03} L. Gurvits, ``What is the finiteness conjecture for linear continuous time
inclusions?'' {\it Proceedings of 42nd IEEE Conference on Decision
and Control (CDC '03),} 2003.

\bibitem{GS07} L. Gurvits, A. Olshevsky, ``Stability Testing of Matrix
Polytopes,'' {\it Proceedings of the European Control Conf.}, July
2007.



\bibitem{HAPL04} D. Henrion, D. Arzelier, D. Peaucelle, J. B. Lasserre, ``On
parameter-dependent Lyapunov functions for robust stability of
linear systems,'' {\it Proceedings of the IEEE Conference on
Decision and Control,} December 2004.


\bibitem{K72} R.M. Karp, ``Reducibility among Combinatorial
Problems,'' In {\it Complexity in Computer Computations,} Ed. R.E.
Miller, J.W. Thatcher, pp. 85-103, Plenum Press, 1972.

\bibitem{M99} V. Monov, ``On the Spectrum of Convex Sets of Matrices,'' {\it IEEE Transactions on Automatic Control,}
vol 44, No. 5, pp. 1009-1012, May 1999

\bibitem{MS65} T.S. Motzkin, E.G. Straus, ``Maxima for graphs and a new proof of a theorem of Turan,'' {\it Canadian
Journal of Mathematics,} Vol. 17, 533-540, 1965.

\bibitem{N93} A. Nemirovskii, ``Several NP-hard problems arising in robust stability
analysis,'' {\it Mathematics of Control, Signals, and Systems,} vol.
6, No. 2, pp. 99-105, 1993.

\bibitem{OP05} R.C.L.F.Oliveira and P.L.D.Peres, ``LMI conditions for the
existence of polynomially parameter-dependent Lyapunov functions
assuring robust stability,'' {\it Proceedings of the 44th IEEE
Conferenceon Decision and Control—European Control Conference,}
December 2005.

\bibitem{PG06} B.T. Polyak, E.N. Gryazina, ``Stability domain in the parameter
space: D-decomposition revisited,'' {\it Automatica,} 42, pp. 13-26,
2006.

\bibitem{PR93} S. Poljak, J. Rohn, ``Checking robust nonsingularity is
NP-hard,'' {\it Mathematics of Control, Signals, and Systems,} Vol.
6, No. 1, pp. 1-9, 1993.

\bibitem{S90} C.B. Soh, ``Schur stability of convex combinations of
matrices,'' {\it Linear Algebra and Its Appl.,} vol. 128, pp.
159-168, 1990.

\bibitem{S05} C.W.Scherer, ``Relaxations for robust linear matrix
inequality problems with verifications for exactness,'' {\it SIAM
Journal on Matrix Analysis and Applications,} 27(2):365–395, 2005.

\bibitem{S86} A. Schrijver, {\it Theory of Linear and Integer
Programming,} John Wiley, 1986.

\bibitem{TB04} P. Tsiotras and P.-A. Bliman, ``An exact stability analysis test
for single-parameter polynomially-dependent linear systems,'' {\it
Proceedings of the 43rd IEEE Conference on Decision and Control,}
December 2004.

\bibitem{val} L.G. Valiant, ``Completeness classes in algebra,''  Proc. of the Eleventh Annual ACM Symposium on Theory of
Computing (Atlanta, Ga., 1979),  pp. 249--261, ACM, New York, 1979.

\bibitem{W91} Q. G. Wang, ``Necessary and sufficient conditions for
stability of a matrix politope with normal vertex matrices,'' {\it
Automatica,} vol. 27, No. 5, pp. 887-888, 1991.



\end{thebibliography}
\end{document}